\documentclass[12pt]{amsart}
\usepackage{latexsym,amsmath,amscd}
\allowdisplaybreaks
\newtheorem{theorem}{Theorem}

\newtheorem{corollary}{Corollary}
\newtheorem{lemma}{Lemma}
\theoremstyle{remark}

\title{Swan Conductors and Torsion in the Logarithmic De Rham Complex}
\date{July 27, 2001}
\author{s\.{i}nan \"{u}nver }
\address{Department of Mathematics\\ University of California-Berkeley\\Berkeley, CA 94720}
\email{sinan@math.berkeley.edu}
\begin{document}
\maketitle
\noindent
\section{Introduction}
In the following $A$ is a discrete valuation ring with maximal ideal 
$\mathfrak{m},$ perfect residue field $k,$ and field of fractions $K,$  
$S=\hbox{Spec}A,$ with closed point $s,$    
and $X/S$ an arithmetic surface over $S$, i.e. an integral, regular 
scheme which 
is proper, flat and of relative dimension one over $S$. We also 
assume that the reduced special fiber
$X_{s,\tiny{\hbox{red}}}$ is a strict normal crossings divisor in $X,$ by this 
we mean that $X_{s,\tiny{\hbox{red}}}$ is a normal crossings divisor in 
$X,$ and 
that the irreducible components of $X_{s,\tiny{\hbox{red}}}$ are regular 
schemes. 

We endow $S$ with the 
log structure corresponding to the natural inclusion 
\begin{eqnarray*}
 \mathcal{O}_{S}- \{ 0 \} \to \mathcal{O}_{S}. 
\end{eqnarray*} 
Similarly, we endow $X$ with the log structure corresponding to the 
natural map 
\begin{eqnarray*}
\mathcal{O}_{X} \cap j_{*}\mathcal{O}^{*}_{X_{K}} \to \mathcal{O}_{X},
\end{eqnarray*} 
where $j:X_{K} \to X$ is the inclusion.
Then the structure map from $X$ to $S$ becomes a 
map of fine log schemes. Let 
$\Omega \dot{ }_{X/S,\tiny{\hbox{log}}}$ denote its logarithmic 
de Rham complex, 
\begin{eqnarray*}
\mathcal{O}_{X} \to \Omega^{1}_{X/S,\tiny{\hbox{log}}} \to \Omega^{2}_
{X/S,\tiny{\hbox{log}}}.
\end{eqnarray*}
Taking $A$-torsion in 
$\Omega \dot{ }_{X/S,\tiny{\hbox{log}}},$ 
we obtain the complex $\Omega \dot{ }_{X/S,\tiny{\hbox{log}},
\tiny{\hbox{tors}}},$ which is a complex supported on the special 
fiber of $X.$ 
We will be mainly interested in the Euler characteristic 
$\chi(\Omega \dot{ }_{X/S,\tiny{\hbox{log}},\tiny{\hbox{tors}}})$ of 
$\Omega\dot{ }_{X/S,\tiny{\hbox{log}},\tiny{\hbox{tors}}},$ which will 
be defined as follows.
If $K \dot{ }$ is a bounded complex of coherent sheaves on $X$ that is exact  
on the generic fiber of $X,$ the hypercohomology groups 
$\hbox{H}\dot{ }(X,K\dot{ })$ are modules of finite length over $A,$ and we 
put 
\begin{eqnarray*}
\chi(K \dot{ })=\sum_{i}(-1)^{i}\hbox{length}_{A}(\hbox{H}^{i}
(X,K \dot{ })).
\end{eqnarray*} 

On the other hand, $X_{K}/K$ has a Swan conductor $\hbox{Sw}(X_{K}/K),$ 
defined as 
follows. Let $K'$ be the strict henselization of the completion of $K$ 
(with respect
to its discrete valuation), and $\ell$ be a prime different from the 
characteristic $p$ 
of $k.$ The action of the wild inertia group of $\hbox{Gal}(\overline{K}'/K')$ 
on $\hbox{H}_{\tiny{\hbox{\'{e}t}}}^{i}(X_{\overline{K}'},\mathbb{Q}_{\ell})$ 
factors through the wild 
inertia group of $\hbox{Gal}(L/K'),$ for a finite Galois extension 
$L$ of $K',$ 
being the continuous action of a pro-$p$  group on a finite dimensional 
$\ell$-adic vector space. 
Let $\hbox{SW}_{L/K'}$ be the Swan module over $\mathbb{Z}_{\ell}$, which is 
the 
unique projective $\mathbb{Z}_{\ell}[\hbox{Gal}(L/K')]$ module having as 
character 
the Swan character of $\hbox{Gal}(L/K').$ Then $\hbox{Sw}_{\overline{K}/K}
(\hbox{H}_{\tiny{\hbox{\'{e}t}}}^{i}(X_{K},\mathbb{Q}_{\ell}))$ is defined 
to be 
\begin{eqnarray*}
\hbox{dim}_{\mathbb{Q}_{\ell}}\hbox{Hom}_{\mathbb{Z}_{\ell}
[\tiny{\hbox{G}}_{L}]}(\hbox{SW}_{L/K'}
,\hbox{H}_{\tiny{\hbox{\'{e}t}}}^{i}(X_{\overline{K}'},\mathbb{Q}_{\ell})), 
\end{eqnarray*}
where $\hbox{G}_{L}=\hbox{Gal}(L/K').$ 
This is independent of the 
choice of $L$ since for a finite Galois extension $L'$ of $K'$ 
containing $L,$  
we have 
\begin{eqnarray*}
\hbox{SW}_{L/K'}=\hbox{SW}_{L'/K'}\otimes_{\mathbb{Z}_{\ell}[\tiny{\hbox{G}}
_{L'}]}
\mathbb{Z}_{\ell}[\hbox{G}_{L}].
\end{eqnarray*}
Finally we define 
\begin{eqnarray*}
\hbox{Sw}(X_{K}/K)=\sum_{i=0}^{2}(-1)^{i}
\hbox{Sw}_{\overline{K}/K}(\hbox{H}_
{\tiny{\hbox{\'{e}t}}}^{i}(X_{K},\mathbb{Q}_{\ell})).
\end{eqnarray*}
Then we have the following theorem. 
\begin{theorem}
With the notation above we have 
\begin{eqnarray*}
\chi(\Omega \dot{ }_{X/S,\tiny{\rm{log}},\tiny{\rm{tors}}})=-\rm{Sw}
\it {(X_{K}/K)}.
\end{eqnarray*}
\end{theorem} 
This can be viewed as the logarithmic version of Bloch's theorem ([Bloch], 
Theorem.1). In fact we 
follow the method in [Bloch] closely, together with the recent work 
of Kato and Saito ([K-S]). 

\section{Proof}
We need the Lemma 1.17 of [K-S] on the local structure of $X/S.$ 
\begin{lemma}
For each point $x \in X_{s}$ there is a Zariski open neighborhood $U$ of $x$ 
in $X$, a 
scheme $P$, etale over $\mathbb{A}^{2}_{S}={\rm Spec} A[t,s]$, and a 
closed immersion 
$i:U \to P$, such that $U$ is defined in $P$ by an equation of the form 
$\pi-ut^{a}s^{b}$, 
where $\pi$ is a uniformizer of $A$, $u$ is a unit in $P$, $t$ and $s$ 
form a system 
of parameters in $\mathcal{O}_{X,x}$ when restricted to $X$, and $a,b$ are 
nonnegative integers.    
\end{lemma}
Proof.(following [K-S])
Since $X_{s,\tiny{\hbox{red}}}$ is a strict normal crossings divisor in $X,$ 
there 
are at most two components of $X_{s,\tiny{\hbox{red}}}$ passing through 
$x.$ If there 
is only one component we choose $t$ to be a local defining function for that 
component, and choose $s$ so that $\{ t,s \}$ is a system of parameters 
at $x.$ If there are two components we choose $t$ and $s$ to be the local 
defining functions for these two components. Now by considering the 
multiplicity of $X_{s}$ along the components of $X_{s,\tiny{\hbox{red}}}$ 
passing  through $x$ we see that $X_{s}$ is defined by $t^{a}s^{b}=0$ in a 
neighborhood of $x$ in $X,$ for some nonnegative integers $a$ and $b.$   
Using $t$ and $s$ we get a map from an open neighborhood $U$ of $x$ to 
$\mathbb{A}^{2}.$ Since $\{s,t\}$ is a system of parameters at $t,$ and the 
residue field $k,$ is perfect, by restricting $U$ if necessary, we may assume 
that this  map is unramified.  
By restricting $U$ if necessary, we can factor this map as a closed immersion 
into a scheme $P$ followed by an etale map from $P$ to $\mathbb{A}^{2}_{S}$ 
([EGA-4], Corollaire 18.4.7).    
Now since $\pi/(t^{a}s^{b})$ is a unit in $U$ by restricting $U$ and $P$ if 
necessary, we can find a unit $u$ in $P$ such that $\pi-ut^{a}s^{b}$ vanishes 
on $U.$ Since $U$ is a divisor in $P,$ and $\pi-ut^{a}s^{b}$ vanishes on 
$U$ to see that $U$ is defined by $\pi-ut^{a}s^{b}$ it suffices to note that 
$\pi-ut^{a}s^{b}$ is not in $\mathfrak{m}^{2}_{P,x}$ (note that since $P$ 
is smooth over $S,$ $\{\pi,t,s \}$ is a system of parameters for 
$\mathcal{O}_{P,x}$). 
\hfill  $\Box$

Continuing with the notation of the lemma and denoting the conormal sheaves 
with $N$ we get an exact sequence, 
\begin{eqnarray*}
0 \to N_{U/P} \to \Omega^{1}_{P/S}|_{U} \to \Omega^{1}_{U/S} \to 0. 
\end{eqnarray*} 
In fact, this sequence is exact for any closed imbedding of $U$ into a 
scheme $P$ smooth over $S.$ 
To see the injectivity of the map $N_{U/P} \to \Omega^{1}_{P/S}|_{U}$ we 
proceed as follows. If we denote the kernel of this map by $M,$ we see that 
since $U_{K}/K$ is smooth, the map is injective over $U_{K},$ and hence 
$M|_{U_{K}}=0.$ 
On the other hand, since the imbedding $U \to P$ is regular $N_{U/P}$ is 
locally free. And therefore $M,$ being a subsheaf of a locally free coherent 
sheaf on an integral scheme $U$ that is supported on a proper closed 
subscheme of $U,$ is zero.  
In the following we endow $P$ with the log structure associated to the 
inclusion 
\begin{displaymath}
\mathcal{O}^{*}_{V} \to \mathcal{O}_{P},
\end{displaymath} 
where $V=P-\{ t^{a}s^{b}=0 \}.$ 
We denote by $\Omega^{1}_{P/S,\tiny{\hbox{log}}}$ the sheaf 
of log differentials where we endow $S$ with the trivial log structure. 
Then for the logarithmic differentials we get a similar exact sequence,
\begin{eqnarray*}
0 \to N_{U/P}\otimes_{A}\mathfrak{m}^{-1} \to 
\Omega^{1}_{P/S,\tiny{\hbox{log}}}|_{U} 
\to 
\Omega^{1}_{U/S,\tiny{\hbox{log}}} \to 0,
\end{eqnarray*} 
where the map 
\begin{displaymath}
\delta:
N_{U/P} \otimes _{A} \mathfrak{m}^{-1} \to \Omega^{1}_{P/S, \tiny{\hbox
{log}}}|_{U}
\end{displaymath}
is the map sending 
\begin{displaymath}
(\pi-ut^{a}s^{b}) \otimes \pi^{-1} \text{ to }
u^{-1} \hbox{d}u + a .\hbox{dlog}t + b .\hbox{dlog}s. 
\end{displaymath}
This resolution of $\Omega^{1}_{U/S,\tiny{\hbox{log}}}$ gives a map 
\begin{eqnarray*}
\alpha_{U}:\Omega^{1}_{U/S,\tiny{\hbox{log}}}
\to \hbox{det}\Omega^{1}_{U/S,\tiny{\hbox{log}}} 
\cong \hbox{Hom}(N_{U/P}\otimes_{A}\mathfrak{m}^{-1},
\Lambda^{2}
\Omega^{1}_{P/S,\tiny{\hbox{log}}}|_{U}),
\end{eqnarray*} 
by the formula 
\begin{eqnarray*}
\alpha_{U} (a)(b)=\widetilde{a} \wedge \delta (b),
\end{eqnarray*}
where $\widetilde{a}$ is a section of $\Omega^{1}_{P/S, \tiny{\hbox{log}}}$ 
that maps to $a.$ 
Since any two resolutions 
of $\Omega^{1}_{U/S,\tiny{\hbox{log}}}$ 
are homotopic and the maps $\alpha_{U}$ for different resolutions are 
compatible with the isomorphisms induced by the homotopies on $\hbox{det}
\Omega^{1}_{U/S,\tiny{\hbox{log}}}$, 
we get a map 
\begin{eqnarray*}
\alpha:\Omega^{1}_{X/S,\tiny{\hbox{log}}} \to \hbox{det}
\Omega^{1}_{X/S,\tiny{\hbox{log}}}.
\end{eqnarray*} 
Let $Z_{U}$ be the closed subscheme of $U$ defined by the 
section of the locally free sheaf $\Omega^{1}_{P/S,\tiny{\hbox{log}}}$ 
corresponding to $\delta.$ As above 
this does not depend on the imbedding, and hence defines a closed subscheme 
$Z$ of $X.$ 
Note 
that $\Omega^{1}_{X/S,\tiny{\hbox{log}}}$ 
is an invertible sheaf over $X-Z,$ and hence $\alpha|_{X-Z}$ is an 
isomorphism. 
Let $C \dot{ }$ 
denote the complex 
\begin{eqnarray*}
\alpha:\Omega^{1}_{X/S,\tiny{\hbox{log}}} \to \hbox{det}
\Omega^{1}_{X/S,\tiny{\hbox{log}}}, 
\end{eqnarray*}
with 
$\Omega^{1}_{X/S,\tiny{\hbox{log}}}$ in degree 1. 

For any bounded complex $K\dot{ }$ of locally free coherent sheaves 
on $X$ which is exact outside a proper closed subscheme $Y$ we have 
a bivariant class 
$$
{\rm ch}^{X}_{Y}(K\dot{ }) \text{ in } 
{\rm A}(Y \to X)_{\mathbb{Q}} 
$$ 
([Fulton], Chapter.18), which is 
the same for any other bounded complex of locally free coherent sheaves  
that is quasi-isomorphic to $K\dot{ }$ 
and exact outside $Y.$ Therefore we can define ${\rm ch}_{Y}^{X}(F)$ 
for any coherent sheaf $F$ supported on $Y.$
If $Y$ is $X_{s}$ we will use the notation 
${\rm ch}_{s}$ for ${\rm ch}^{X}_{X_{s}}.$
Now by the Riemann-Roch theorem ([Fulton], Chapter.18; [Saito], Lemma.2.4) 
we have 
$$
\chi(\Omega \dot{ }_{X/S,\tiny{\hbox{log,tors}}})=
{\rm deg}({\rm ch}_{s}(\Omega \dot{ }_{X/S,\tiny{\hbox{log,tors}}})
\cap {\rm Td}({X/S})),
$$
where ${\rm Td}(X/S)$ is the Todd class of $X/S.$ Since $C\dot{ }$ is 
exact outside $Z,$ after choosing any resolution of $\Omega^{1}_{X/S
,\tiny{\hbox {log}}}$ we can define ${\rm ch}_{s}(C\dot{ }).$  
Using the following lemma, we will work with $C \dot{ }$ instead of 
$\Omega \dot{ }_{X/S,\tiny{\hbox{log}},\tiny{\hbox{tors}}}.$ 
\begin{lemma}
With the notation above, we have 
\begin{eqnarray*}
{\rm ch}_{s}(C \dot{ })={\rm ch}_{s}(\Omega \dot{ }
_{X/S,{\rm log},{\rm tors}}).
\end{eqnarray*} 
\end{lemma}
Proof. Note that $\hbox{ker}(\alpha)=\Omega^{1}_{X/S,\tiny{\hbox{log}},
\tiny{\hbox{tors}}},$ since 
$\hbox{det}\Omega^{1}_{X/S,\tiny{\hbox{log}}}$ 
is an invertible sheaf and $\alpha$ is an isomorphism 
on the generic fiber. Therefore to finish the proof of the lemma we need to 
show that 
$\hbox{ch}_{s}(\hbox{coker}(\alpha))=\hbox{ch}_{s}
(\Omega^{2}_{X/S,\tiny{\hbox{log,tors}}}).$ First 
note that 
\begin{displaymath}
\hbox{coker}(\alpha) \cong \hbox{det}\Omega^{1}_{X/S,
\tiny{\hbox{log}}}\otimes \mathcal{O}_{Z}.
\end{displaymath} 
To see this we can work locally and choose an imbedding of $U$ as in Lemma.1. 
Let 
\begin{displaymath}
u^{-1}\hbox{d}u+a.\hbox{dlog}t+b.\hbox{dlog}s=(a+t.x)
\hbox{dlog}t+(b+s.y)\hbox{dlog}s,
\end{displaymath} 
for some $x,y$ in $\mathcal{O}_{U}.$ Then 
$Z$ is defined in  $U$ by the ideal 
\begin{eqnarray*}
&&(a+t.x,b+s.y) \text{, if } a \neq 0 \text{ and } b \neq 0 \text{, or by }\\ 
&&(a+t.x,y) \text{, if } a \neq 0 \text{ and } b=0.
\end{eqnarray*}
$\Omega^{1}_{X/S,\tiny{\hbox{log}}}$ is generated by  ${\rm d}t,$ ${\rm d}s,$ 
$\hbox{dlog}t,$ (if $a \neq 0$), and $\hbox{dlog}s$ (if $b \neq 0$) 
subject to the relations 
\begin{eqnarray*}
&&(a+t.x)\hbox{dlog}t+(b+s.y)\hbox{dlog}s=0,\\ 
&& t.{\rm dlog}t={\rm d}t, \text{ and }s.{\rm dlog}s={\rm d}s.
\end{eqnarray*}
Note that we may view $\hbox{det}\Omega^{1}_{X/S,\tiny{\hbox{log}}}$ as a 
subsheaf of $\Omega^{1}_{X_{K}/K,\tiny{\hbox{log}}}.$ 
If $a \neq 0$ and $b \neq 0$  ${\rm det} \Omega^{1}_{X/S,\tiny{\hbox{log}}}$
 is generated by 
\begin{eqnarray*}
&& \frac{1}{b+s.y}\hbox{dlog}t  \text{, if } b+s.y \neq 0 
    \text{, or by }\\ 
&& \frac{1}{a+t.x}\hbox{dlog}s  \text{, if } 
   a+t.x \neq 0  . \\ 
\end{eqnarray*}
 Assume without loss of generality 
that $b+s.y$ is nonzero. Then the image of $\Omega^{1}_{X/S,
\tiny{\hbox{log}}}$ in 
$\hbox{det}\Omega^{1}_{X/S,\tiny{\hbox{log}}}$ 
is generated by 
\begin{displaymath}
\hbox{dlog}t=(b+s.y)\frac{1}
{b+s.y}\hbox{dlog}t,\text{ and } \hbox{dlog}s=(a+t.x)\frac{1}{b+s.y}
\hbox{dlog}t.
\end{displaymath} 
Therefore the cokernel of $\alpha$ 
is $\hbox{det}\Omega^{1}_{X/S,\tiny{\hbox{log}}}\otimes \mathcal{O}_{Z}.$ 
If $a \neq 0$ and $b=0$ then 
${\rm det}\Omega^{1}_{X/S,\tiny{\hbox{log}}}$ is generated by 
\begin{eqnarray*}
&& \frac{1}{y}{\rm dlog}t, \text{ if } y \neq 0, \text{ or by}\\
&& \frac{1}{a+t.x}{\rm dlog}s \text{, if } a+t.x \neq 0, 
\end{eqnarray*} 
in $\Omega^{1}_{X_{K}/K,\tiny{\hbox{log}}},$ 
and  we  similarly arrive at the conclusion. 

Next, if $a \neq 0$ and $b \neq 0,$ $\Omega^{2}_{X/S,\tiny{\hbox{log}}}$  
is generated by $\hbox{dlog}t \wedge 
\hbox{dlog}s$ with the relations 
\begin{displaymath}
(a+t.x)\hbox{dlog}t \wedge \hbox{dlog}s=0 \text{, and } 
(b+s.y) \hbox{dlog}t \wedge \hbox{dlog}s=0. 
\end{displaymath}
And if $a \neq 0$ and $b=0,$ $\Omega^{2}_{X/S,\tiny{\hbox{log}}}$ is 
generated by ${\rm dlog}t \wedge {\rm d}s$ with the relations 
$$
(a+t.x){\rm dlog}t \wedge {\rm d}s=0, \text{ and } y{\rm dlog}t \wedge 
{\rm d}s=0.  
$$
This shows that $\Omega^{2}_{X/S,\tiny{\hbox{log}}}$ is an 
invertible sheaf on $Z.$ And using again the local desription we see that 
$\Omega^{1}_{X/S,\tiny{\hbox{log}}}|Z$ 
is locally free of rank 2, and we have 
\begin{displaymath}
\Omega^{2}_{X/S,\tiny{\hbox{log,tors}}}
=\Omega^{2}_{X/S,\tiny{\hbox{log}}}=\Lambda^{2}
\Omega^{1}_{X/S,\tiny{\hbox{log}}}|_{Z}.
\end{displaymath} 
Restricting the resolution of $\Omega ^{1}_{X/S,\tiny{\hbox{log}}}$ 
over $U$ to $Z_{U}$, we obtain 
\begin{eqnarray*}
0 \to L^{1}i^{*}\Omega^{1}_
{X/S,\tiny{\hbox{log}}}|_{Z_{U}} \to N_{U/P}|_{Z_{U}} \otimes_{A} 
\mathfrak{m}^{-1} \to \Omega^{1}_{P/S,\tiny{\hbox{log}}}|_{Z_{U}} 
\to \Omega^{1}_{U/S,\tiny{\hbox{log}}}|_{Z_{U}} \to 0, 
\end{eqnarray*} 
where $i:Z \to X$ is the inclusion. Here the second and the fourth arrows 
are isomorphisms. In particular, we have 
$$
L^{1} i^{*}\Omega^{1}_{X/S,\tiny{\hbox{log}}}|_{Z_{U}} 
\cong N_{U/P}|_{Z_{U}} \otimes_{A}\mathfrak{m}^{-1}.  
$$
The exact sequence also shows that 
\begin{eqnarray*}
\Lambda^{2}\Omega^{1}_{X/S,\tiny{\hbox{log}}}|_{Z} 
\cong \hbox{det}\Omega^{1}_{X/S,\tiny{\hbox{log}}}|_{Z}\otimes
L^{1}i^{*}\Omega^{1}_{X/S,\tiny{\hbox{log}}}.
\end{eqnarray*}
Using a filtration of $\mathcal{O}_{Z}$ with graded pieces supported on 
integral subschemes of $Z,$ 
we see that  
to prove the lemma it is enough to show that $L^{1}i^{*}
\Omega^{1}_{X/S,\tiny{\hbox{log}}}|_{T}\cong \mathcal{O}_{T},$ for any 
integral 
curve $T$ in $Z.$ For the rest of the 
proof we use the method of the proof of Proposition 3.1 in [Saito], 
in this very 
explicit (and easier)  case.   
First note that if 
$k:U \to Q$ is a closed immersion with $Q$ smooth over $S,$ then $k$ is a 
regular 
immersion. Since the inclusion $T_{U} \to U$ is also a regular immersion, 
we have an 
exact sequence of locally free sheaves on $T_{U},$
\begin{eqnarray*}
0 \to N_{U/Q}|_{T_{U}} \to N_{T_{U}/Q} \to N_{T_{U}/U} \to 0.
\end{eqnarray*}     
And similarly we have an exact sequence,
\begin{eqnarray*}
0 \to N_{Q_{s}/Q}|_{T_{U}} \to N_{T_{U}/Q} \to N_{T_{U}/Q_{s}} \to 0,
\end{eqnarray*} 
in particular $N_{Q_{s}/Q}|_{T_{U}} \to N_{T_{U}/Q}$  is injective. 
Furthermore 
for an immersion $U \to P$ as in Lemma.1 we claim that 
\begin{eqnarray*}
0 \to N_{P_{s}/P}|_{T_{U}} \to N_{T_{U}/P} \to N_{T_{U}/U} \to 0
\end{eqnarray*}
is exact. To see this we only need to check that 
\begin{displaymath}
N_{P_{s}/P}|_{T_{U}}={\rm ker} (N_{T_{U}/P} \to N_{T_{U}/U}).
\end{displaymath} 
Since this is a local question on $X,$ by restricting 
$U$ and $P$ we will assume that $T_{U}$ is defined by $t$ on $U$. 
Denoting $\pi-ut^{a}s^{b}$ by $g,$ we need to show that 
\begin{displaymath}
\pi/(\pi^{2},\pi t)={\rm ker}((t,g)/(t,g)^{2} \to (t,g)/(t^{2},g)).
\end{displaymath}
Since by assumption $T$ is contained in 
$Z,$ $X/S$ 
is not smooth along $T$, and so $a \geq 2.$ Therefore $ut^{a}s^{b}
 \in (t^{2}),$ and  
 $\pi \in (t^{2},g).$ This shows that $\pi/(\pi^{2},\pi t)$ is in the kernel. 
To see the converse we only need to note that, 
\begin{displaymath}
g-\pi=ut^{a}s^{b}=0 \text{ in } (t,g)/(t,g)^{2},
\end{displaymath} 
since $ut^{a}s^{b} \in (t^{2}).$ 
This proves the claim. Tensoring the exact sequence with $\mathfrak{m}^{-1}$ 
and observing that 
$$
N_{P_{s}/P}|_{T_{U}}\otimes_{A}\mathfrak{m}^{-1}
\cong \mathcal{O}_{T_{U}}, 
$$
we obtain the exact sequence 
$$
0 \to \mathcal{O}_{T_{U}} \to N_{T_{U}/P} \otimes_{A}\mathfrak{m}^{-1} 
\to N_{T_{U}/U}\otimes_{A}\mathfrak{m}^{-1} \to 0.  
$$
On the other hand, using the isomorphism 
$$
L^{1} i^{*}\Omega^{1}_{X/S,\tiny{\hbox{log}}}|_{T_{U}}
\cong N_{U/P}|_{T_{U}} \otimes_{A}\mathfrak{m}^{-1}
$$
and the exact sequence 
$$
0 \to N_{U/P}|_{T_{U}}\otimes_{A}\mathfrak{m}^{-1} \to 
N_{T_{U}/P}\otimes_{A}\mathfrak{m}^{-1} \to N_{T_{U}/U}\otimes_{A}
\mathfrak{m}^{-1} \to 0, 
$$
we get an exact sequence 
$$
0 \to L^{1} i^{*}\Omega^{1}_{X/S,\tiny{\hbox{log}}}|_{T_{U}} 
\to N_{T_{U}/P} \otimes_{A}\mathfrak{m}^{-1} \to N_{T_{U}/U}
\otimes_{A}\mathfrak{m}^{-1} \to 0. 
$$
Therefore we see that 
\begin{displaymath}
L^{1}i^{*}\Omega^{1}_{X/S,\tiny{\hbox{log}}}|_{T_{U}}\cong 
\mathcal{O}_{T_{U}} 
\end{displaymath}
by viewing them both as the kernel of 
\begin{displaymath}
N_{T_{U}/P}\otimes_{A}\mathfrak{m}^{-1} 
\to N_{T_{U}/U}\otimes_{A}\mathfrak{m}^{-1}. 
\end{displaymath}
If we take imbeddings of $U$ into $P$ and $P'$ as above, taking 
$Q=P\times P'$ we get the 
inclusions
\begin{eqnarray*}
&&L^{1}i^{*}\Omega^{1}
_{X/S,\tiny{\hbox{log}}}|_{T_{U}}\to N_{U/Q}|_{T_{U}}\otimes_{A}
\mathfrak{m}^{-1}
\to N_{T_{U}/Q}\otimes_{A}\mathfrak{m}^{-1},\\ 
&&\text{ and } \mathcal{O}_{T_{U}}\cong  N_{Q_{s}/Q}|_{T_{U}} \otimes_{A} 
\mathfrak{m}^{-1}
\to N_{T_{U}/Q}\otimes_{A}\mathfrak{m}^{-1}.
\end{eqnarray*}  
Using this we see that the isomorphism 
$L^{1}i^{*}
\Omega^{1}_{X/S,\tiny{\hbox{log}}}|T_{U}
\cong \mathcal{O}_{T_{U}}$ does not depend on the choice of the local 
imbedding satisfying 
Lemma.1. Therefore we obtain 
\begin{displaymath}
L^{1}i^{*}\Omega^{1}_{X/S,\tiny{\hbox{log}}}|_{T}\cong 
\mathcal{O}_{T}.
\end{displaymath}
This finishes the proof of the lemma.
\hfill $\Box$  

Using Lemma.2 we see that 
\begin{eqnarray*}
\chi(\Omega \dot{ }_{X/S,\tiny{\hbox{log,tors}}})=
\hbox{deg}(\hbox{ch}_{s}(\Omega \dot{ }
_{X/S,\tiny{\hbox{log,tors}}}) \cap \hbox{Td}(X/S))=
\hbox{deg}(\hbox{ch}_{s}(C \dot{ })\cap 
\hbox{Td}(X/S)).
\end{eqnarray*}
Let 
\begin{displaymath}
F\dot{ }:0 \to F_{m}\to \cdots \to F_{0}\to 0 
\end{displaymath}
be a complex of locally 
free coherent sheaves on $X,$ which is exact outside a proper closed set 
$Y.$ This exact 
sequence on $X-Y$ gives a canonical trivialization over $X-Y,$ of the line 
bundle 
\begin{displaymath}
\hbox{det}F\dot{ }=\otimes_{0 \leq i \leq m}(\hbox{det}F_{i})^
{\otimes(-1)^{i}}.
\end{displaymath} 
This gives a rational section $s$ of $\hbox{det}F\dot{ }.$ 
Denote  the image  of the 
divisor of $s$ in the  Chow group $\hbox{A}_{*}Y$ of $Y,$ by 
$\gamma.$ We  will need the following lemma.
\begin{lemma}
We have the equality
\begin{eqnarray*}
{\rm ch}_{Y,1}^{X} (F\dot{ })\cap [X]=\gamma, \text{ in } {\rm A}_{*}Y.
\end{eqnarray*}
\end{lemma}    
Proof. Let $f_{i}=\hbox{rank}F_{i},$ $F_{-1}=0,$ and 
\begin{displaymath}
G_{i}=\hbox{Grass}_{f_{i}}(F_{i}\oplus
F_{i-1}), \text{ the Grasmannian of } f_{i}  
\text{ planes in } F_{i}\oplus F_{i-1},
\end{displaymath} 
 for $0 \leq i \leq m.$ 
 Let $\xi_{i}$ be the tautological subbundle, of rank $f_{i},$ of $F_{i}\oplus 
F_{i-1}$ on $G_{i}.$ 
Let 
\begin{displaymath}
G=G_{m} \times \cdots \times G_{0}, \text { with the projections }
p_{i}:G \to G_{i}, \text{ and } \pi:G \to X.
\end{displaymath}
Let 
\begin{displaymath}
\xi=\sum_{i=0}^{m}(-1)^{i}p_{i}^{*}\xi_{i}, \text{ and } 
\hbox{det}\xi=\otimes_{0 \leq i \leq m}(\hbox{det}p_{i}^{*}\xi_{i})^
{\otimes(-1)^{i}}.
\end{displaymath}
Furthermore we denote the kernel of $d_{i}:F_{i}\to F_{i-1}$ by $K_{i},$ 
of rank $k_{i},$ 
and $H_{i}=\hbox{Grass}_{k_{i}}(F_{i}).$ 
Finally, let $W$ be the closure of $\varphi(X\times
\mathbb{A}^{1})$ in $G\times \mathbb{P}^{1},$ 
where $\varphi:X\times \mathbb{A}^{1}
\to G\times \mathbb {A}^{1}$ is the map 
sending $(x,\lambda)$ to 
\begin{displaymath}
(\hbox{Graph}(\lambda \cdot d_{m}(x)),\cdots,
\hbox{Graph}
(\lambda \cdot d_{0}(x)),\lambda ).
\end{displaymath} 
Over $\varphi 
((X-Y)\times \mathbb{A}^{1}),$ $\hbox{det}\xi$ has a natural trivialization 
since 
\begin{displaymath}
\hbox{det}p_{i}^{*}\xi_{i} \cong \hbox{det}\pi^{*}  K_{i}\otimes 
\hbox{det} \pi^{*} K_{i-1} 
\text{ over } \varphi((X-Y)\times \mathbb{A}^{1}). 
\end{displaymath}
This trivialization gives a divisor, say $D,$ on $\varphi 
(X\times \mathbb{A}^{1}).$  
$D$ is 
supported on $G_{Y}\times \mathbb{A}^{1}.$ 
\begin{displaymath}
\pi_{*}([D_{0}])=\gamma \text{ in } \hbox{A}_{*}Y 
\end{displaymath}
Let $t \in \mathbb{P}^{1}-\{0,\infty \}.$ As 
$$
\pi_{*}([D_{0}])=\pi_{*}([D_{t}]),
$$ 
we will be done if we can show that $\pi_{*}([D_{t}])$ is equal to 
$\hbox{ch}_{Y,1}^{X}(F\dot{ })\cap [X].$ For this, by definition, we need to 
show that  
\begin{displaymath}
\pi_{*}([D_{t}])=\pi_{*}(\hbox{ch}_{1}(\xi)\cap [T]), 
\end{displaymath}
where $[T]=[W_{\infty}]-[\widetilde{X}],$ 
and $\widetilde{X}$ is the irreducible component of $W_{\infty}$ which projects
 birationally onto 
$X.$ Note that over $X-Y,$ $\varphi$ can be extended to a function 
 $$
 \varphi:(X-Y)\times \mathbb{P}^{1} \to G \times \mathbb{P}^{1}
 $$ 
as follows. 
For $(x,[\lambda_{0},\lambda_{1}]) \in (X-Y) \times \mathbb{P}^{1},$ and 
$0 \leq i \leq m,$   
let $\varphi_{i}'(x,[\lambda_{0},
\lambda_{1}])$ denote the point of $G_{i}$ over $x$ that corresponds to 
 $$
 \{(v_{i},v_{i-1}) \in F_{i}(x)\oplus K_{i-1}(x):\lambda_{0}v_{i-1}=
\lambda_{1}d_{i}v_{i}\} \subseteq  F_{i}(x)\oplus F_{i-1}(x).
$$ 
Then 
$$
\varphi(x,[\lambda_{0},\lambda_{1}])=(\varphi_{m}'(x,[\lambda_{0},
\lambda_{1}]),\cdots,\varphi_{0}'(x,[\lambda_{0},\lambda_{1}]),
[\lambda_{0},\lambda_{1}]).
$$ 
And $\widetilde{X}$ is the 
closure of $\varphi((X-Y)\times \{\infty\})$ in $G\times \{\infty \}.$ 
Now as 
$$
\hbox{det}p_{i}^{*}\xi_{i}\cong \hbox{det}\pi^{*}K_{i} \otimes 
\hbox{det}\pi^{*}
K_{i-1} \text{ over } \varphi((X-Y) \times (\mathbb{P}^{1}-\{0\})),
$$  
$\hbox{det}\xi$ has a natural trivialization over $\varphi((X-Y)
\times (\mathbb{P}^{1}-\{ 0 \})),$ which gives a divisor, say $D',$  
on $W-W_{0},$ supported on $G_{Y}\times (\mathbb{P}^{1}-\{0\}).$ If 
$t \in \mathbb{P}^{1}-\{0,\infty\}$ then $[D_{t}]=[D_{t}']$ in 
$\hbox{A}_{*}W_{t}$ being divisors associated to the 
same line bundle $\hbox{det}\xi|_{W_{t}}.$ Noting that 
$[D'_{t}]=[D'_{\infty}],$ we only need to show that 
$$       
\pi_{*}(D'_{\infty})=\pi_{*}(D'.T) \text{ in } 
\hbox{A}_{*}Y,
$$
or that 
$$
\pi_{*}(D_{\infty}'. \widetilde{X} )=0 \text{ in } \hbox{A}_{*}Y.
$$
If 
\begin{displaymath}
\Psi :X-Y \to H \text{ is the map that sends } 
x  \text{ to } (K_{m}(x),\cdots,K_{0}(x)),  \text{ and} 
\end{displaymath}
\begin{displaymath}
\iota: H \to G \text{ is the map that sends } 
(V_{m},\cdots,V_{0}) \text{ to } (V_{m}\oplus V_{m-1},
\cdots, V_{0}),
\end{displaymath}
then 
\begin{displaymath}
\widetilde{X}=\iota(\overline{\Psi(X-Y)}).
\end{displaymath}
However, as $\iota^{*}(D'_{\infty})$ is the divisor 
corresponding to a section of 
\begin{displaymath}
\otimes_{0 \leq i \leq m}\hbox{det}\iota^{*}
\xi_{i} ^{\otimes (-1) ^{i}}
\cong \otimes_{0 \leq i \leq m}(\hbox{det}q_{i}^{*}\zeta_{i}\otimes 
\hbox{det} q_{i-1}^{*}\zeta_{i-1})^{\otimes (-1)^{i}} \cong 
\mathcal{O}_{H}, 
\end{displaymath}
that is nonzero on $H_{X-Y},$ 
where $\zeta_{i}$ is the tautological subbundle of $H_{i},$ and 
$$q_{i}:
H \to H_{i}$$ 
is the projection, we see that $$\pi_{*}(D'_{\infty}.
\widetilde{X})=0.$$ 
This finishes the proof of the lemma.  
\hfill $\Box$  

Let $0 \to E_{m} \to \cdots \to E_{1} \to \Omega^{1}_{X/S,\tiny{\hbox{log}}} 
\to 0$ be a 
resolution 
of $\Omega^{1}_{X/S,\tiny{\hbox{log}}}$ by locally free 
sheaves of finite rank. Now 
consider the complex 
\begin{eqnarray*}
E \dot{ }:0 \to E_{m}\to \cdots \to E_{1} \to \hbox{det}\Omega^{1}_{X/S,
\tiny{\hbox{log}}} \to 0,
\end{eqnarray*}
where the map $E_{1} \to \hbox{det}\Omega^{1}_{X/S,
\tiny{\hbox{log}}}$ is the composition of 
the differential 
\begin{displaymath}
d_{1}:E_{1} \to \Omega^{1}_{X/S,\tiny{\hbox{log}}}, 
\end{displaymath}
and the canonical map 
\begin{displaymath}
\alpha:\Omega^{1}_{X/S,\tiny{\hbox{log}}}\to \hbox{det}\Omega^{1}_{X/S,
\tiny{\hbox{log}}}. 
\end{displaymath}
Applying Lemma.3 to $E\dot{ }$ we obtain
\begin{eqnarray*}
\chi(\Omega\dot{ }_{X/S,\tiny{\hbox{log,tors}}})&=&
\hbox{deg}(\hbox{ch}_{s}(C\dot { })\cap \hbox{Td}(X/S))\\
&=&\hbox{deg}(\hbox{ch}_{s}(E\dot{ })\cap \hbox{Td}(X/S))=
\hbox{deg}(\hbox{ch}_{s,2}(E\dot{ }) \cap [X]).
\end{eqnarray*} 
We will need the following lemma. 
\begin{lemma}
We have ${\rm ch}_{s,2}{(E \dot{ })\cap [X]}={\rm c}_{s,2}
(\Omega^{1}_
{X/S,\tiny{{\rm log}}}) 
\cap [X]$ in $({\rm A}_{*}X_{s})_{\mathbb{Q}}.$
\end{lemma} 
 Proof. First of all note that we have $\hbox{ch}_{s,1}(E \dot{ }) \cap [X]=0$ 
 by Lemma.1, 
hence 
$$
\hbox{ch}_{s,2}(E \dot{ }) \cap [X]=\frac{\hbox{c}_{s,1}^{2}(E\dot{ })\cap [X]}
{2} - \hbox{c}_{s,2}(E\dot{ })\cap [X]=-\hbox{c}_{s,2}(E \dot{ })\cap [X].
$$ 
Let $E'$ 
denote the complex $0 \to E_{m} \to \cdots \to E_{1} \to 0,$ where we put 
$E_{1}$ in 
degree 0. We use the notation in Lemma.3, where the objects with $'$ refer to 
those 
corresponding to $E'.$ 

We have 
\begin{eqnarray*}
\hbox{ch}_{s,2}(E \dot{ }) \cap [X] &=& -\hbox{c}_{s,2}(E \dot{ })\cap [X]=
\pi_{*}(-\hbox{c}_{2}(\sum_{i=0}^{m}(-1)^{i}[p_{i}^{*}\xi_{i}]) \cap [T]) 
\text{, and}\\
&=& \pi_{*}((-\hbox{c}_{2}(\sum_{i=1}^{m}(-1)^{i}[p_{i}^{*}\xi_{i}])-
\hbox{c}_{1}
(\sum_{i=1}^{m}(-1)^{i}[p_{i}^{*}\xi_{i}])\cdot 
\hbox{c}_{1}(p_{0}^{*}\xi_{0})\\
&&-\hbox{c}_{2}(p_{0}^{*}\xi_{0}))\cap [T]).
\end{eqnarray*}  
Since $\xi_{0}$ is a line bundle this is equal to 
\begin{eqnarray*}
\pi_{*}((-\hbox{c}_{2}(\sum_{i=1}^{m}(-1)^{i}[p_{i}^{*}\xi_{i}])&-&
\hbox{c}_{1}(\sum_{i=1}^{m}(-1)^{i}[p_{i}^{*}\xi_{i}])\cdot \hbox{c}_{1}(p_{0}
^{*}\xi_{0}))\cap [T])\\
&=& \pi_{*}((-\hbox{c}_{2}(\sum_{i=1}^{m}(-1)^{i}[p_{i}^{*}\xi_{i}])+
\hbox{c}_{1}^{2}(p_{0}^{*}\xi_{0}))\cap [T])\\
&&-\pi_{*}((\hbox{c}_{1}(\sum_{i=0}^{m}(-1)^{i}[p_{i}^{*}\xi_{i}])\cdot 
\hbox{c}_{1}(p_{0}^{*}\xi_{0}))\cap [T])\\
&=& \pi_{*}((-\hbox{c}_{2}(\sum_{i=1}^{m}(-1)^{i}[p_{i}^{*}\xi_{i}])+
\hbox{c}_{1}^{2}(p_{0}^{*}\xi_{0}))\cap [T])\\
&& -(\hbox{c}_{s,1}(E \dot{ })\cap [X])\cdot \hbox{c}_{1}(\xi_{0}) \quad 
\text{(since $p_{0}=\pi $).}
 \end{eqnarray*}
Lemma.3 shows that this is equal to 
\begin{eqnarray*}
&& \pi_{*}((-\hbox{c}_{2}(\sum_{i=1}^{m}(-1)^{i}[p_{i}^{*}\xi_{i}])+
\hbox{c}_{1}^{2}(p_{0}^{*}\xi_{0}))\cap [T])\\
&=& \pi_{*}((-\hbox{c}_{1}^{2}(\sum_{i=1}^{m}(-1)^{i+1}[p_{i}^{*}\xi_{i}])
+\hbox{c}_{2}(\sum_{i=1}^{m}(-1)^{i+1}[p_{i}^{*}\xi_{i}])\\
&&+\hbox{c}_{1}^{2}(p_{0}^{*}\xi_{0}))\cap [T]).
\end{eqnarray*}
If $f:G \to G'$ is the projection, note that we have $f^{*}p_{i}^{'*}\xi_{i}'
\cong p_{i}^{*}\xi_{i},$ for $2 \leq i \leq m.$ Therefore the last expression 
is equal to 
\begin{eqnarray*}
&&\pi_{*}((-\hbox{c}_{1}^{2}(\sum_{i=1}^{m}(-1)^{i+1}[p_{i}^{*}\xi_{i}])+ 
\hbox{c}_{2}(\sum_{i=1}^{m}(-1)^{i+1}[f^{*}p_{i}^{'*}\xi_{i}'])\\
&+& c_{1}(\sum_{i=1}^{m}(-1)^{i+1}[f^{*}p_{i}^{'*}\xi_{i}'])\cdot \hbox{c}_{1}
([p_{1}^{*}\xi_{1}]-[f^{*}p_{1}^{'*}\xi_{1}'])\\
&&+\hbox{c}_{2}([p_{1}^{*}\xi_{1}]
-[f^{*}p_{1}^{'*}\xi_{1}'])+\hbox{c}_{1}^{2}(p_{0}^{*}\xi_{0}))\cap [T]). 
\end{eqnarray*}
This is equal to 
\begin{eqnarray*}
\pi_{*}(-\hbox{c}_{1}^{2}(p_{0}^{*}\xi_{0})&+& c_{1}([p_{0}^{*}\xi_{0}]-
[p_{1}^{*}\xi_{1}]+[f^{*}p_{1}^{'*}\xi_{1}'])\cdot \hbox{c}_{1}([p_{1}^{*}
\xi_{1}]-[f^{*}p_{1}^{'*}\xi_{1}'])\\
&+& \hbox{c}_{2}([p_{1}^{*}\xi_{1}]-[f^{*}p_{1}^{'*}\xi_{1}'])+c_{1}^{2}(
p_{0}^{*}\xi_{0}))\cap [T]) +\hbox{c}_{s,2}(\Omega^{1}_{X/S,
\tiny{\hbox{log}}})\cap [X]\\
&=& \pi_{*}((\hbox{c}_{1}([\pi^{*}\hbox{det}\Omega^{1}_{X/S,
\tiny{\hbox{log}}}]-[p_{1}^{*}
\xi_{1}]+[\pi^{*}E_{1}])\cdot \hbox{c}_{1}([p_{1}^{*}\xi_{1}]-[\pi^{*}E_{1}])\\
&& +\hbox{c}_{2}([p_{1}^{*}\xi_{1}]-[\pi^{*}E_{1}]))\cap [T]) + \hbox{c}_{s,2}
(\Omega^{1}_{X/S,
\tiny{\hbox{log}}})\cap [X].\\
&& \text{(cf. the proof of Lemma.3)}
\end{eqnarray*}
Now, since $p_{1}^{*}\xi_{1}$ is the tautological subbundle of $\pi^{*}E_{1}
\oplus \pi^{*}\hbox{det}\Omega^{1}_{X/S,\tiny{\hbox{log}}}$ on $G,$ we have an 
exact sequence 
$ 0 \to p_{1}^{*}\xi_{1} \to \pi^{*}E_{1}\oplus \pi^{*}\hbox{det}\Omega^{1}
_{X/S,\tiny{\hbox{log}}} \to Q \to 0,$ with $Q$ a line bundle. Therefore 
we obtain 
\begin{eqnarray*}
\hbox{c}_{2}([p_{1}^{*}\xi_{1}]-[\pi^{*}E_{1}])&=& \hbox{c}_{2}([
\pi^{*}\hbox{det}\Omega^{1}_{X/S,\tiny{\hbox{log}}}]-[Q])\\
&=& \hbox{c}_{2}(\pi^{*}\hbox{det}\Omega^{1}_{X/S,\tiny{\hbox{log}}})-
\hbox{c}_{1}
([\pi^{*}\hbox{det}\Omega^{1}_{X/S,\tiny{\hbox{log}}}]-[Q])\cdot 
\hbox{c}_{1}(Q)- \hbox{c}_{2}(Q).
\end{eqnarray*} 
Since $\xi_{0}$ and $Q$ are line bundles  
\begin{displaymath}
\hbox{c}_{1}(Q)=
\hbox{c}_{1}([\pi^{*}E_{1}]+[\pi^{*}\hbox{det}\Omega^{1}_{X/S,
\tiny{\hbox{log}}}]-[p_{1}^{*}\xi_{1}]), \text{ and}
\end{displaymath} 
\begin{eqnarray*}
\hbox{c}_{2}([p_{1}^{*}\xi_{1}]-[\pi^{*}E_{1}])=-\hbox{c}_{1}([p_{1}^{*}
\xi_{1}]-[\pi^{*}E_{1}])\cdot \hbox{c}_{1}([\pi^{*}E_{1}]+[\pi^{*}\hbox{det}
\Omega^{1}_{X/S,\tiny{\hbox{log}}}]-[p_{1}^{*}\xi_{1}]). 
\end{eqnarray*}  
Combining this with the expression for $\hbox{ch}_{s,2}(E\dot{ })\cap [X]$ 
above we obtain 
\begin{eqnarray*}
\hbox{ch}_{s,2}(E\dot { })\cap [X]=\hbox{c}_{s,2}(\Omega^{1}_
{X/S,\tiny{\hbox{log}}})\cap [X].
\end{eqnarray*}
\hfill $\Box$

Using Lemma.4 we obtain 
\begin{eqnarray*}
\chi(\Omega \dot{ }_{X/S,\tiny{\hbox{log,tors}}})=\hbox{deg}
(\hbox{ch}_{s,2}(E\dot{ }) \cap
 [X])= \hbox{deg}(\hbox{c}_{s,2}(\Omega^{1}_{X/S,\tiny{\hbox{log}}})\cap [X]).
\end{eqnarray*}
Now the logarithmic version of Bloch's conductor formula  ([K-S], Theorem.1.15) 
says that 
\begin{eqnarray*}
\hbox{deg}(\hbox{c}_{s,2}(\Omega^{1}_{X/S,
\tiny{\hbox{log}}})\cap[X])=-\hbox{Sw}(X_{K}/K).
\end{eqnarray*}
Combining this  with the above we obtain
\begin{eqnarray*}
\chi(\Omega\dot{ }_{X/S,\tiny{\hbox{log,tors}}})=-\hbox{Sw}(X_{K}/K).
\end{eqnarray*} 
\hfill $\Box$

We now give a consequence of the proof of Theorem.1. In the following  
if $C$ is a 0-dimensional subscheme of $X,$ and $[G]$ and $[H]$ are 
curves in $X,$  
we denote by ${\rm deg}[C]$ the degree of $C$ with respect to $k,$ 
and by $[G]\cdot [H]$ the intersection number of the curves $G$ and $H.$ 
Let 
$D \subseteq X$ be a curve supported on the special fiber $X_{s},$ $L$ a 
line bundle on $X,$ and $K$ and $E$ divisors on $X$ such that 
$\mathcal{O}(K)\cong \hbox{det}\Omega^{1}_{X/S},$ the dualizing sheaf 
of $X/S,$ and $\mathcal{O}(E) \cong L.$ Then we have the following lemma.
\begin{lemma}
We have the equality 
\begin{eqnarray*}
\chi(i_{*}(L|_{D}))=[E]\cdot [D] -\frac{1}{2}\chi(\omega_{D/k}),
\end{eqnarray*}
where $i:D \to X$ is the inclusion.
 \end{lemma}
Proof. By the Riemann-Roch theorem we have 
\begin{eqnarray*}
\chi(i_{*}(L|_{D}))=\hbox{deg}(\hbox{ch}_{D}^{X}(L|_{D}) \cap  
\hbox{Td}(X/S)).
\end{eqnarray*}
By applying the Riemann-Roch theorem to the closed immersion 
$i:D \to X$ we obtain 
\begin{eqnarray*}
\hbox{ch}^{X}_{D}(i_{*}(L|_{D}))&=&i_{*}(\hbox{ch}_{D}^{D}(L|_{D})
\cdot  \hbox{Td}(N_{D/X}\check{ })^{-1})\\
&=&(1+\hbox{c}_{1}(L|_{D}))\cdot(1+\frac{1}{2}\hbox{c}_{1}(\mathcal{O}(-D)|_{D}
)).
\end{eqnarray*}
Combining this with the above we obtain 
\begin{eqnarray*}
\chi(i_{*}(L|_{D}))=[E]\cdot [D] -\frac{([K]+[D])\cdot [D]}{2}.
\end{eqnarray*}
Finally using the adjunction formula for $D \to X$ we obtain the expression 
in the statement of the lemma.
\hfill $\Box$

Let $Z_{1}$ and $Z_{2}$ denote the closed subschemes of $Z$ consisting of 
the components of $Z$ which have codimension 1 and codimension 2 in
$X$ respectively. With this notation, we have the following corollary. 
\begin{corollary}
We have the following equality
\begin{eqnarray*}
{\rm Sw}(X_{K}/K)=\chi(\Omega^{1}_{X/S,{\rm log},{\rm tors}})+
(2[Z_{1}]-[Z_{1,{\rm red}}])\cdot
[Z_{1}]-\frac{1}{2}\chi(\omega_{Z_{1}/k})-{\rm deg}[Z_{2}].
\end{eqnarray*}
\end{corollary}
Proof. Using Theorem.1 we see that we only need to prove the equality 
\begin{eqnarray*}
\chi(\Omega^{2}_{X/S,\tiny{\hbox{log}}})=([Z_{1,\tiny{\hbox{red}}}]
-2[Z_{1}])\cdot [Z_{1}]+
\frac{1}{2}\chi(\omega_{Z_{1}/k})+{\rm deg}[Z_{2}].
\end{eqnarray*}
 The proof of Lemma.2 shows that 
\begin{eqnarray*}
\chi(\Omega^{2}_{X/S,\tiny{\hbox{log}}})=\chi(\hbox{det}
\Omega^{1}_{X/S,\tiny{\hbox{log}}}|_{Z}).
\end{eqnarray*} 
Then we have 
\begin{eqnarray*}
\chi(\Omega^{2}_{X/S,\tiny{\hbox{log}}})
=\chi(\hbox{det}\Omega^{1}_{X/S,\tiny{\hbox{log}}}|_{Z_{1}})
+\hbox{deg}[Z_{2}].
\end{eqnarray*}
Using the lemma above we obtain that  
\begin{eqnarray*}
\chi(\Omega^{2}_{X/S,\tiny{\hbox{log}}})=
[K_{\tiny{\hbox{log}}}]\cdot [Z_{1}]- \frac{1}{2}\chi(\omega_{Z_{1}/k})
+{\rm deg}[Z_{2}],
\end{eqnarray*}
where $K_{\tiny{\hbox{log}}}$ is a divisor on $X$ such that $\mathcal{O}
(K_{\tiny{\hbox{log}}}) \cong {\rm det}\Omega^{1}_{X/S,\tiny{\hbox{log}}}.$  
Now we also have  
\begin{eqnarray*}
[K_{\tiny{\hbox{log}}}]=[K]-[Z_{1}]+[Z_{1,\tiny{\hbox{red}}}]. 
\end{eqnarray*} 
To see this we only need to look at the multiplicities at codimension 
1 points. Using the notation in the proof of Lemma.2, viewing 
${\rm det}\Omega^{1}_{X/S}$ as a subsheaf of 
$$
\Omega^{1}_{X_{K}/K}\cong \Omega^{1}_{X_{K}/K,\tiny{\hbox{log}}},
$$ 
it is generated by 
\begin{eqnarray*}
&&\frac{1}{t^{a-1}y}{\rm dlog}t, \text{ if } y \neq 0, \text{ or by }\\
&&\frac{1}{t^{a-1}(a+t.x)}{\rm d}s, \text{ if } a+t.x \neq 0,
\end{eqnarray*}
in a Zariski neighborhood of a point with $a \neq 0$ and $b=0.$ 
Using the similar description of ${\rm det}\Omega^{1}
_{X/S,\tiny{\hbox{log}}}$ in the proof of Lemma.2, we arrive at the 
formula as claimed above.   
Using this and the adjunction formula we obtain that 
\begin{eqnarray*}
\chi(\Omega^{2}_{X/S,\tiny{\hbox{log}}})=
([Z_{1,\tiny{\hbox{red}}}]-2[Z_{1}])\cdot [Z_{1}]+
\frac{1}{2}\chi(\omega_{Z_{1}/k})+{\rm deg}[Z_{2}].
\end{eqnarray*} 
\hfill $\Box$

\end{document}